\documentstyle[twoside]{article}

\newtheorem{thm}{\indent{\sc Theorem}}[section]
\newtheorem{defn}[thm]{\indent{\sc Definition}} 
 
\newtheorem{por}[thm]{\indent{\sc Porism}}
\newtheorem{prop}[thm]{\indent{\sc Proposition}} 
\newtheorem{cor}[thm]{\indent{\sc Corollary}}
\newtheorem{ex}{\indent {\sc Example}}[section]

\newcommand{\ssect}{\subsection}
\newcommand{\sssect}{\subsubsection}

\newcommand{\eqref}[1]{equation~(\ref{#1})}
\newcommand{\eref}[1]{\eqref{#1}}
\newcommand{\Eqref}[1]{Equation~(\ref{#1})}
\newcommand{\Eref}[1]{\Eqref{#1}}
\newcommand{~}{\nolinebreak[3] }

\newcommand{\seper}{\nopagebreak \begin{center}
	\underline{\hspace{2in}}
	\end{center}}

\newcommand{\bold}[1]{{\em (#1)}\index{#1}}

 \newcommand{\proof}[1]{\proo{#1}$\Box $}
\newcommand{\proo}[1]{{\em Proof: }#1}
\renewcommand{\Box}{ \hspace{1cm}{\rule{1.2ex}{2ex}}}

\newcommand{\rn}[1]{\left\lfloor #1 \right\rceil}

\newcommand{\phil}{\nonumber \\*[-1mm]  &&\: }

\newcommand{\SB}[1]{\mbox{{\bf \scriptsize #1}}}

\newcommand{\digress}[1]{ \begin{quotation}
	 #1 {\em End of Digression.} \end{quotation}}

\newcommand{\fig}[2]{\begin{table}[htbp]
	\caption{#1} 
	\begin{center}\fbox{ \scriptsize \begin{tabular}#2 \end{tabular}}
	\end{center}
\end{table}}

\begin{document}
\include{title}

\begin{abstract}
We  generalize the Stirling numbers of the first kind $s(a,k)$ to the case
where $a$ may be an arbitrary real number. In particular, we study the
case in which $a$ is an integer. There, we discover new combinatorial
properties held by the classical Stirling numbers, and analogous properties
held by the Stirling numbers $s(n,k)$ with $n$ a negative integer.
\seper
{\bf G\'{e}n\'{e}ralisation des nombres de Stirling}

On g\'{e}n\'{e}ralise ici les nombres de Stirling du premier ordre
$s(a,k)$ au cas o\`u $a$ est un r\'eel quelconque. On s'interesse en
particulier au cas o\`u $a$ est entier. Ceci permet  de mettre en
evidence de nouvelles propri\'et\'es combinatoires aux quelles
obeissent les nombres de Stirling usuels et des propri\'et\'es
analougues auquelles obeissent les nombres de Stirling $s(n,k)$ o\`u
$n$ est un entier n\`egatif.
\end{abstract}

\begin{center}
{\it Dedicated to\\
Robert and Julia}
\end{center}

\tableofcontents

\section{Introduction}

The Stirling numbers are some of the most important combinatorial constants
known. It is the hope of this paper, to generalize the Stirling numbers of the
first kind $s(n,k)$ to the case where $n$ need not be a nonnegative integer.
This leads to \eref{stun} which points out the link between the two types
of Stirling numbers.

When $n$ is an integer, we discover many  new combinatorial
properties held by the classical Stirling numbers, and analogous properties
held by the Stirling numbers $s(n,k)$ with $n$ a negative integer.
However, we  defer some of their interesting combinatorial properties to
another paper \cite{hybrid}.

Finally, one should refer to \cite{ch4} for one of the most important
applications of these new constants: the 
calculation of the harmonic logarithms which form a basis for the Iterated
Logarithmic Algebra. 

\ssect{The Lower Factorial}\label{DisEx}

Recall that for $n$ a nonnegative integer, one defined the lower factorial
function $(x)_{n}$ to be the product
$$ (x)_{n} = x(x-1)\cdots (x-n+1).$$
Similarly, for $n$ a negative integer,  we can define 
$$(x)_n = \prod _{i=k}^{-1}(y-i)^{-1} =
{1}/{(y+1)(y+2)\cdots (y-k)}$$

Equivalently, $(y)_{k}$ can be defined recursively by the
requirements:
\begin{eqnarray*}
(y)_{0}&=&1\\
(y)_{k}&=&(y-k+1)(y)_{k-1}\mbox{ for all $k$.}
\end{eqnarray*}

\digress{We make a brief {\em  digression} on notation in order to prevent
any possible confusion. We use the symbol $(x)_{n}$ to
denote the ``falling powers'' $y^{\underline{n}}=
y(y-1)\cdots (y-n+1)$; however, 
many researchers---for example, Askey\index{Askey{,} Richard} and
Henrici\index{Henrici}---reserve
this notation for the ``rising powers'' $y^{\overline{n}} =
y(y+1)\cdots (y+n-1)$. Actually, as Knuth\index{Knuth} has pointed out,
Pochhammer\index{Pochhammer} who 
devised this notation did not intend either of these
definitions; he used $(y)_{n}$ to denote $y(y-1)\cdots (y-n+1)/n!$.}

For $a$ a real number, we could define $(y)_{a}$ in terms of the Gamma
function. 

\begin{defn}\bold{Lower Factorial}\label{lowerC} 
Let $a$ be a real number, and
let $y$ be a formal variable. The define the lower factorial by
$$ (y)_{a}=\frac{\Gamma (y+1)}{\Gamma (y-a+1)} $$
where $\Gamma (y)$ denotes the formal power series expansion in the variable
$y$ of the  Gamma function.\index{Gamma Function}
\end{defn}

Note that for all real numbers $a$, 
\begin{equation}\label{Det1}
(y)_{a}=(y-a+1)(y)_{a-1}.
\end{equation}

\digress{We {\em digess} to indicate how one could proceed more formally.
Require \cite{Ueno} for all real
numbers $a$ that
$$ \sum _{i\in \SB{Z}}\frac{(y)_{a+i}}{\rn{a+i}!}t^{a+i} =
\sum _{i\in \SB{Z}}\frac{y^{a+i}}{\rn{a+i}!}\log
(1+t)^{a+i}$$
where $y$ and $t$ are {\em variables.} It then
follows (see \cite{Char} or \cite{UP}) that
$$ (y)_{a+i} = y^{a+i}(1+A) $$
where 
$$ A=\sum_{k\geq 1}y^{-k} \sum _{j=1}^{k}
(-1)^{j}{a+i+j-1\choose j+k }{a+i+k\choose k-j }S(j+k,j)$$
and the $S(j+k,j)$ are the {\em Stirling numbers of the second
kind}\index{Stirling Numbers of the Second Kind}. 
Note that $A$ is not a Laurent series so the
inner summation does not give the Stirling numbers of the
first kind.}

\ssect{The Stirling Numbers of the First Kind}

We  generalize the classical definition of the Stirling numbers of
the first kind, and  derive some remarkable identities satisfied by them.

\begin{defn}
\bold{Stirling Numbers of the First Kind}  For all real numbers 
$a$ and for all nonnegative integers $k$, we define the
{\em Stirling number of the first kind} $s(a,k)$ of degree $a$ and order $k$
$s(a,k)$ to be the coefficient $[y^{k}](y)_{a}$ in the Taylor expansion of the
lower factorial. Thus,
$$ (y)_{a}=\sum _{k\geq 0}s(a,k)y^{k}. $$
\end{defn}

{Note that for $a$ a 
positive integer, this corresponds to the usual definition of
Stirling numbers of the first kind.}

\begin{ex}\begin{enumerate}
\item For $n$ a nonnegative integer, $s(n,k)$ is the usual Stirling number of
the first kind. That is, $(-1)^{n+k}s(n,k)$ is the number of permutations of
$n$ elements which is the product of $k$ disjoint cycles. 
\item For all nonnegative integers $k$, $s(0,k)=\delta _{0,k}.$
\item $s(a,0)=1$ except when $a$ is a positive integer in which case
$s(a,0)=0$.
\item See Table~\ref{harmfig}. 
\end{enumerate}
\end{ex}
\fig{Stirling Numbers of the First Kind, $s(n,k)$}{{c|rrrrr|rrrrrrr}
\label{harmfig}
\normalsize $k \backslash n$& \normalsize $-5$ &
\normalsize $-4$ & \normalsize $-3$ & \normalsize $-2$ & \normalsize $-1$ &
\normalsize 0 & \normalsize 1 & \normalsize 2 &
\normalsize 3 & \normalsize 4 & \normalsize 5 & \normalsize $-6$ \\
\hline
\normalsize 0& $\frac{1}{120}$ & $\frac{1}{24}$ & $\frac{1}{6}$ & $\frac{1}{2}$
& 1 & 1 & 0&0&0&0&0&0\\[0.1in]
\normalsize 1& $-\frac{137}{7200}$ & $-\frac{25}{288}$ & $-\frac{11}{36}$ &
$-\frac{3}{4}$ & $-1$ & 0 & $1$ & $-1$ & $2$ & $-6$ & $24$ & $-120$ \\[0.1in]
\normalsize 2& $\frac{12,019}{432,000}$ & $\frac{415}{3456}$ & $\frac{85}{216}$
& $\frac{7}{8}$ & 1 & 0&0 & 1 & $-3$ & 11 & $-50$ & 274 \\[0.1in]
\normalsize 3& $-\frac{874,853}{25,920,000}$ & $-\frac{5845}{41,472}$ &
$-\frac{575}{1296}$ & $-\frac{15}{16}$ & $-1$ & 0&0&0 & 1 & $-6$
 & 35 & $-225$\\[0.1in]
\normalsize 4& $\frac{58,067,611}{1,555,200,000}$ & $\frac{76,111}{497,664}$ &
$\frac{3661}{7776}$ & $\frac{31}{32}$ & 1 & 0&0&0&0& 1 & $-10$ & 85\\[0.1in]
\normalsize 5& $-\frac{3,673,451,957}{93,312,000,000}$ &
$-\frac{952,525}{5,971,968}$ & $-\frac{22,631}{46,656}$ & $-\frac{63}{64}$ &
$-1$ & 0&0&0&0&0& 1 & $-15$\\[0.1in]
\normalsize 6& $\frac{226,576,031,859}{5,598,720,000,000}$&
$\frac{11,679,655}{71,663,616}$ & $\frac{137,845}{279,936}$ & $\frac{127}{128}$
&1&0&0&0&0&0&0&1}

\begin{thm}\label{lemon}
For all $a$ and for all positive integers $k$,
\begin{equation}\label{LemEq}
s(a+1,k)=s(a,k-1)-as(a,k).
\end{equation}
\end{thm}

\proof{\Eref{Det1}.}
\section{Properties of the Stirling Numbers of the First Kind}
\sssect{Nonpositive Degree}
\index{Stirling Numbers,Nonpositive Degree}\label{nonnegsec} 

In this section, we derive several identities which hold for
Stirling numbers of the first kind with degrees which are nonpositive integers.

Recall that a {\em linear partition}\index{Partition{, }Linear}
$\rho $ is a nonincreasing infinite
sequence, $(\rho _{i})_{i\geq 1}$,  of nonnegative integers
which is eventually zero. 
For example, $\rho =(17,2,2,1,0,0\ldots )$ is a linear partition.
Each nonzero $\rho_{i}$ is called a {\em part}\index{Parts} $\rho$. In the
above example,  the multiset of parts of $\rho$ is $\{1,2,2,17\}$. The number
of parts of $\rho$ is denoted $\ell (\rho )$.
A linear partition $\rho$ is said to be a {\em partition of $n$}
if the sum of its parts is $n$, and we write $\rho \vdash n$.  
The product of the parts of $\rho $ is denoted by $\pi
(\rho)$. 

The set of all linear partitions is denoted by ${\cal P}$. A
linear partition is said to have {\em distinct parts}\index{Parts,Distinct} if
its multiset of parts is, in fact, a set. The set of all linear
partitions with distinct parts is denoted by ${\cal P}^{*}.$

The zero sequence is a linear partition. It has no
parts. Thus, it is a partition of zero with distinct parts,
and the product of its parts is one.

As promised, Proposition~\ref{olddefn} and its corollaries give
enumerative interpretations of the Stirling numbers
of nonnegative degree. 

\begin{prop}\label{olddefn}
\bold{Harmonic Relation}\label{HarmRel}
For any nonnegative integers $n$ and $k$, the 
Stirling number of degree $n$ and order $k$ is given by the sum
\begin{eqnarray*}
s(-n,k)&=&\frac{(-1)^{k}}{n!}\displaystyle \sum_{\scriptstyle \rho \in {\cal
P}\atop {\scriptstyle \ell (\rho)=k\atop\scriptstyle \rho
_{1}\leq n } }\pi(\rho)^{-1};
\end{eqnarray*}
over all linear partitions 
$\rho $ with $k$ parts and with
no part greater than $n$.
\end{prop}

We offer two proofs. First,

{\em Proof 1:} Let $n$ and $k $ be as above, and define 
$$d_{n}^{(k)}=\displaystyle \sum_{\scriptstyle \rho \in {\cal
P}\atop {\scriptstyle \ell (\rho)=k\atop\scriptstyle \rho
_{1}\leq n } }\pi(\rho)^{-1}.$$ 
By Theorem~\ref{lemon}, it suffices to verify the recursion
\begin{equation}\label{lem}
nd_{n}^{(k)}-d_{n}^{(k-1)}=nd_{n-1}^{(k)}
\end{equation}
for $n$ and $k$ positive,
since the boundary conditions are easy to verify.

Consider the following series of equalities.
\begin{eqnarray*} 
nd_{n}^{(k)}-d_{n}^{(k-1)}&=&
\left(\displaystyle n \sum _{\scriptstyle \mu \in {\cal P}
\atop {\scriptstyle \ell (\mu
)=k \atop \scriptstyle\mu _{1} \leq n} } \displaystyle \pi(\mu
)^{-1}\right)
 - \left(\sum _{\scriptstyle \mu \in {\cal P} \atop
{\scriptstyle \ell (\mu )=k-1
\atop \scriptstyle\mu _{1} \leq n} } \displaystyle \pi(\mu
)^{-1}\right)\\*[0.1in]
&=& \left(\displaystyle n \sum _{\scriptstyle \mu \in {\cal
P} \atop {\scriptstyle \ell
(\mu )=k \atop \scriptstyle\mu _{1} \leq n} } \displaystyle \pi(\mu
)^{-1} \right)
- \left(n \sum _{\scriptstyle \mu \in {\cal P} \atop
{\scriptstyle \ell (\mu )=k
\atop \scriptstyle\mu _{1} = n} } \displaystyle \pi(\mu
)^{-1}\right)\\[0.1in] 
&=& n
\displaystyle \sum _{\scriptstyle \mu \in {\cal P} \atop
{\scriptstyle \ell (\mu
)=k\atop \scriptstyle\mu _{1} < n} } \displaystyle \pi(\mu )^{-1}\\*[0.1in]
&=&nd_{n-1}^{(k)}.\end{eqnarray*}
Thus, \eref{lem} holds.$\Box $

Alternately, we could adopt the following proof due to Y. C.
Chen.\index{Chen{,} Y. C.}

{\em Proof 2:} For $n$ positive,
\begin{eqnarray*}
S(-n,k)&=& [y^{k}](y)_{-n}\\*
&=& \frac{1}{n!}[y^{k}]\frac{1}{(1+y)(1+y/2)\cdots (1+y/n)}\\
&=& \frac{(-1)^{k}}{n!}[y^{k}] \left(\sum _{\rho _{1}\geq 0}y^{\rho
_{2}}\right) \left(\sum _{\rho _{2}\geq 0}(y/2)^{\rho
_{1}}\right)\cdots \left(\sum _{\rho _{n}\geq
0}(y/n)^{\rho _{n}}\right)  \\*
&=&\frac{(-1)^{k}}{n!}\displaystyle \sum_{\scriptstyle \rho \in
{\cal 
P}\atop {\scriptstyle \ell (\rho)=k\atop\scriptstyle \rho_{1}\leq n }
}\pi(\rho)^{-1}.\Box  \end{eqnarray*}

From the second proof, we have the following the Porism involving the {\em
complete symmetric function}\index{Complete Symmetric Function}  
$h_{n}(x_{1},x_{2},\ldots ) $. It is defined explicitly by
the sum 
\begin{equation}\label{early h}
 h_{n}(x_{1},x_{2},\ldots ) = \sum _{\scriptstyle 1\leq \alpha _{1}\leq \alpha
_{2}\leq \cdots \leq \alpha _{n}} \prod
_{k=1}^{n}x_{\alpha _{k}},
\end{equation}
and it is defined implicitly by the generating function
$$ \prod _{n\geq 1}(1-x_{n}y)^{-1}=\sum _{n\geq
0}h_{n}(x_{1},x_{2},\ldots )y^{n}. $$

\begin{por}
\label{complete}
Let $n$ and $k$ be nonnegative integers. Then
$$ s(-n,k) = h_{k}(-1,-1/2,\ldots ,-1/n)/n!.\Box  $$
\end{por}

Thus, we see from the Harmonic Relation (Proposition~\ref{HarmRel}) that
the Stirling numbers are simply related to the partial
sums of the harmonic series. For $k=1$,
Proposition~\ref{olddefn} yields a sum over partitions of
length one with no part greater than $n$. There are $n$ such
partitions; they are the integers from $1$ to $n$. Thus, the sum is
the sum of the reciprocals of the first $n$ integers, so
that 
\begin{equation}\label{harmsum}
n!s(-n,1)=-\left(1+\frac{1}{2}+\cdots +\frac{1}{n} \right)
\end{equation}

For the Stirling numbers of order 2, we obtain similarly:
\begin{equation}\label{hsum2}
n!s(-n,2)=1+\frac{1}{2}\left(1+\frac{1}{2}\right)+
\frac{1}{3}\left(1+\frac{1}{2}+\frac{1}{3}\right)+\cdots 
+\frac{1}{n}\left(1+\frac{1}{2}+\cdots +\frac{1}{n}\right),
\end{equation}
and for order 3, we obtain
\begin{eqnarray*}
-n!s(-n,2)&=&1
	\phil +\frac{1}{2} \left[1+ \frac{1}{2}
\left(1+\frac{1}{2}\right)\right]
	\phil +\frac{1}{3}\left[1+ \frac{1}{2}
\left(1+\frac{1}{2}\right) +\frac{1}{3}
\left(1+\frac{1}{2}+\frac{1}{3}\right)\right] 
	\phil \vdots
	\phil +\frac{1}{n}\left[1+ \frac{1}{2}
\left(1+\frac{1}{2}\right) +\frac{1}{3}
\left(1+\frac{1}{2}+\frac{1}{3}\right)+\cdots
+\frac{1}{n}\left(1+\frac{1}{2}+\cdots +\frac{1}{n}\right)\right].
\end{eqnarray*}

Proposition~\ref{helpful}
turns out to be very useful in the calculation of Stirling
numbers of the first kind.

\begin{prop}\label{helpful}
\bold{Knuth} Let $n$ and $k$ be nonnegative integers (not
both zero). Then $s(n,k)$ is given by the 
following finite sum:
\begin{equation}\label{HelpEq}
s(-n,k)= \frac{ (-1)^{k+1}}{n!} \sum_{m=1}^{n} {n\choose m} (-1)^{m} m^{-k}.
\end{equation}
\end{prop}

\proo{By consideration of the examples above, the
Proposition holds for $n=0$ and $k=0$. Now, by induction we need only
show that \eref{lem} holds for the summation on the  right side of
\eref{HelpEq}.
\begin{eqnarray*}
\lefteqn{\sum _{m\geq 1} {n-1\choose m } (-1)^{m}m^{-k}
+\frac{1}{n} \sum_{m\geq 1} {n\choose m } (-1)^{m} m^{1-k}}\\*
 &= & \sum _{m\geq 1} (-1)^{m} m^{-k} \left({n-1\choose m
}+\frac{m}{n}{n\choose m } \right)\\
 &= & \sum _{m\geq 1} (-1)^{m} m^{-k} \left({n-1\choose m
}+{n-1\choose m-1 } \right)\\*
 &= & \sum _{m\geq 1} (-1)^{m} m^{-k}{n-1\choose m}.\Box 
\end{eqnarray*}}

We would like to generalize to
the case where  $n$ need not be a nonnegative integer; however,
this is impossible since in this case the sum is not only {\em infinite} but
{\em divergent}. 

Note that Proposition~\ref{helpful} is the analog of the following classical
result \cite{Char} involving Stirling numbers of the {\em second} kind
$S(k,n)$.
$$ S(k,n)=\frac{(-1)^{n}}{n!}\sum_{m=1}^{k} {n\choose m}(-1)^{m}m^{k}. $$
Thus, in some sense we can say that 
\begin{equation}\label{stun}
S(k,n)=(-1)^{n+k+1}s(-n,-k).
\end{equation}

Proposition~\ref{olddefn} has several more corollaries:

\begin{cor}\label{more}
For nonnegative integers $n$ and $k$, the Stirling number of degree $n$ and
order $k$ is given by the sum
$$ s(-n,k)= \frac{(-1)^{k}n}{(n-1)!}\displaystyle \sum_{\scriptstyle \rho \in
{\cal P}\atop {\scriptstyle \ell (\rho)=k+1\atop\scriptstyle \rho_{1}= n }}
\pi(\rho )^{-1}$$  
 over all linear partitions 
$\rho $ with $k+1$ parts and whose largest part is $n$.
\end{cor}

\proof{Add the part $n$ to each partition $\rho$ being summed
over in Corollary~\ref{olddefn}.}

\begin{cor}\label{Stirl}
Let $n$ and $k$ be nonnegative integers. Then the Stirling number of the first
kind of degree $-n$ and order $k$ is given by the sums
\begin{eqnarray*}
s(-n,k)
	&=&\frac{(-1)^{k}}{n!}\displaystyle \sum _{\scriptstyle M \subseteq 
\{1,2,\ldots,n \}
\atop \scriptstyle |M|=k }
\left(\prod_{m\in M}m^{-1}\right)\\*
	&=&\frac{(-1)^{k}}{(n-1)!}\displaystyle \sum _{\scriptstyle M \subseteq
\{1,2,\ldots ,n \}\atop{\scriptstyle |M|=k+1 \atop \scriptstyle n\in M }}
\left(\prod_{m\in M}m^{-1}\right)
\end{eqnarray*}
 over {\em multisets} $M$ where all of the products
are computed with the proper  multiplicities.
\end{cor}

\proof{Every linear partition is associated with a unique
multiset of positive numbers called its parts. In the identities from 
Proposition~\ref{olddefn} and Corollary~\ref{more}, sum over these
multisets instead of the partitions themselves.}

\begin{cor}
Let $n$ and $k$ be nonnegative integers. Then the Stirling number of the first
kind of degree $-n$ and order $k$ is given by the sum
$$s(-n,k) = \frac{(-1)^{k}}{n!} \displaystyle \sum 
\left(\prod_{i=1}^{n}i^{-m_{i}}\right)$$
 over all sequences
$(m_{i})_{i=1}^{n}$ of $n$ nonnegative integers
which sum to $k$.
\end{cor}

\proof{Every linear partition is determined by the
number $m_{i}$ of times each integer $i$ occurs as a part.
Hence, we can sum over sequences of nonnegative integers
$m_{i}$.} 

We note that
\begin{equation}\label{star39}
\lim_{k \rightarrow +\infty} s(-n,k)= (-1)^{k}/(n-1)!
\end{equation}
for all
nonnegative $n$. 

\sssect{Positive Degree}
\index{Stirling Numbers,Positive Degree}\label{negsec} 

Now, we develop the classical Stirling numbers (those of positive degree)
in a similar vein.

Recall that the trivial partition has no parts, and therefore the
product of its parts is one; however, there are no partitions with
$-1$ parts.

\begin{prop}\label{ptn}\bold{Harmonic Relation}
Let $k$ be a nonnegative integer, and let $n$ be a positive
integer. Then the Stirling number of the first kind of degree $n$ and order
$k$ is given by the sums
\begin{eqnarray*}
s(n,k)& = & (-1)^{n+k}(n-1)! \displaystyle \sum _{\scriptstyle \mu \in
{\cal P}^{*}\atop{\scriptstyle \ell(\mu)=k-1 \atop \scriptstyle \mu
_{1}<n}} \pi(\mu )^{-1}\\*
&=& (-1)^{n+k} \displaystyle \sum _{\scriptstyle \mu \in
{\cal P}^{*}\atop{\scriptstyle \ell(\mu)=n-k \atop \scriptstyle \mu_{1}<n}}
\pi(\mu ),
\end{eqnarray*}
over all linear partitions $\mu$ with 
$k-1$ parts all of which are distinct and less than $n$.
\end{prop}

\proo{Note that
\begin{eqnarray*}
s(n,k) & = & [x^{k}](x)_{n}\\*
&=& [x^{k}]\prod _{i=0}^{n-1}(x-i)\\
&=&[x^{k}] \sum _{\scriptstyle \mu \in {\cal P}^{*} \atop \scriptstyle \mu
_{1}<n } (-1)^{\ell (\mu )}x^{n-\ell (\mu )}\pi (\mu )\\
&=& [x^{k}] \sum _{\scriptstyle \mu \in {\cal
P}^{*} \atop \scriptstyle \mu _{1}<n } 
\frac{(n-\ell (\mu))!}{(n-\ell (\mu )-k)!}
 (-1)^{\ell (\mu )} x^{n-\ell(\mu )-k} \pi (\mu )\\
&=& \sum _{\scriptstyle \mu \in {\cal
P}^{*} \atop {\scriptstyle \mu _{1}<n \atop \scriptstyle
\ell (\mu )=n-k}} k! (-1)^{n-k} \pi (\mu )\\*
&=& (-1)^{n+k}(n-1)!\sum _{\scriptstyle \nu \in {\cal
P}^{*} \atop {\scriptstyle \nu _{1}<n \atop \scriptstyle
\ell (\nu )=k-1}} \pi (\nu )^{-1}.\Box 
\end{eqnarray*}}

Note also that $s(n,k)=0$ if $k>n>0$ or if $k=0$ and $n>0$,
since there is no partition with $k$ distinct parts all less than $k$,
or with $-1$ parts.

By way of example, let
 us consider the extreme cases. If $k=n>0$, then we must
have $k-1$ 
distinct parts less than $k$. There is only one way to do
this; we must 
use the partition consisting of the integers from 1 through
$k-1$. Thus, 
$s(k,k)=1$.

Conversely, for $k=1$ and $n>0$, we sum over partitions with no
parts. The trivial partition is the only such partition, so
$s(n,1)=(-1)^{n-1}(n-1)!$.

The Stirling numbers of positive degree and order 2 are related to the partial
sums of  the harmonic series
\begin{equation}\label{-sum}
s(n,2)=(-1)^{n}(n-1)!\left(1+\frac{1}{2}+\cdots+\frac{1}{n-1}  \right).
\end{equation}
Again, as happened for Stirling numbers of
nonnegative degree, higher orders correspond to
generalizations of the harmonic series.

Some useful equivalent formulations of Proposition~\ref{ptn} follow:
\begin{cor}\label{ptn2}
Let $k $ be a nonnegative integer, and let $n $ be a
positive integer. Then
$$ s(n,k)= (-1)^{n+k}n! \displaystyle -n\sum 
_{\scriptstyle \mu \in {\cal P}^{*}\atop
{\scriptstyle \ell(\mu)=k 
\atop \scriptstyle \mu_{1}=n}}\pi(\mu )^{-1}.$$
\end{cor}

\proof{Add the part $n$ to each partition $\mu $ being summed
over in Proposition~\ref{ptn}.}

\begin{cor}\label{unord}
Let $k$ be a nonnegative integer, and let $n$ be a positive
integer. Then
$$\begin{array}{rclcl}
s(n,k)
&=& \displaystyle (-1)^{n+k}(n-1)! \sum_{\scriptstyle  S \subseteq
\{1,2,\ldots ,n-1 \} \atop \scriptstyle |S|=k-1}
\left( \prod_{s\in S}s^{-1}\right)
&=& \displaystyle \sum_{\scriptstyle  S \subseteq
\{1,2,\ldots ,n-1 \} \atop \scriptstyle |S|=n-k}
\left( \prod_{s\in S}(-s) \right)\\*
&=&\displaystyle (-1)^{n+k}n! \sum 
_{{\scriptstyle S \subseteq \{1,2,\ldots ,n \}
\atop{\scriptstyle |S|=k \atop \scriptstyle n\in S } }}
\left(\prod _{s\in S}s^{-1}\right)
\end{array} $$
where the sums range over sets $S$.
\end{cor}

\proof{These identities can be obtained from
Proposition~\ref{ptn} and Corollary~\ref{ptn2} by summing over the set
of parts of $\mu $ instead of $\mu $ itself.}

\begin{cor}
For all nonnegative integers $k$ and all positive integers $n$, the Stirling
number of the first kind of degree $n$ and order $k$ is given by the sum
$$ s(n,k)= 
\displaystyle (-1)^{n+k}n! \sum _{\scriptstyle \mu \vdash n \atop
\scriptstyle \ell (\mu )=k}
\left(\pi (\mu)^{-1}\right)
\left(\prod _{j=1}^{n}m_{i}(\mu)!\right)^{-1}$$
 over partitions $\rho $ 
of the number $n$ into exactly $k$
parts, and $m_{i}(\rho)$ denote the number of times $i$ occurs as part
of $\rho $.
\end{cor}

\proof{$|s(n,k)|$ is the number of permutations of
$n$ letters with $k$ cycles.
The number of permutations of cycle type $\rho $ is $n!\left(\prod_{i\geq
1}i^{m_{i}}m_{i}! \right)^{-1}$ where $m_{i}$ is the number of parts of $\rho $
equal to $i$.} 

In contrast to \eref{star39},
\begin{prop}\label{2ways}
For $n$ positive, 
\begin{equation}
\sum _{k=1}^{-n}(-1)^{k}s(n,k)=(-1)^{n+1}/(n-1)!.
\end{equation}
Furthermore,
\begin{equation}\label{suffice2}\label{suffice}
\sum _{k\geq 0}(-1)^{k}s(n,k)=(-1)^{n+1}/(n-1)!. 
\end{equation}
\end{prop}

{\em Proof 1:} It  suffices to demonstrate \eref{suffice},
since $s(n,k)= 0$ if $k=0$ or if $k>n$.
Let us expand the left hand side of the
\eref{suffice} as follows:
\begin{eqnarray*}
\sum _{k\geq 0}(-1)^{k}s(n,k) &=&(-1)^{n}\sum _{k\geq 0} |s(-n,k)|\\*
&=&(-1)^{n}\sum _{k\geq 0} \left| \left\{\pi \in
S_{-n}: \pi \mbox{ has $k$ cycles} \right\} \right|\\*
&=&(-1)^{n}n!\Box 
\end{eqnarray*}

Alternately, we have the following proof.

{\em Proof 2:}As mentioned above, it suffices to verify
\eref{suffice2}.
\begin{eqnarray*}
\sum_{k\geq 0}s(n,k)x^{k}&=& (y)_{n}\\*
 \sum_{k\geq 0} s(n,k) (-1)^{k}&=& (-1)(-1-1)\cdots (-1-n+1)\\*
&=& (-1)^{n}n!.\Box 
\end{eqnarray*}

Recall the definition of the {\em elementary symmetric
function}\index{Elementary Symmetric Function}
$e_{n}(x_{1},x_{2},\ldots ) $. It is defined explicitly by
the sum
\begin{equation}\label{early e}
e_{n}(x_{1},x_{2},\ldots ) = \sum _{\scriptstyle 0<\alpha _{1}<\alpha
_{2}<\cdots <\alpha _{n}} \prod_{k=1}^{n} x_{\alpha_{k}}
\end{equation}
 all linear partitions $\mu $ with
distinct parts, 
and it is defined implicitly by the generating function
$$ \prod _{n\geq 1}(1+x_{n}y)=\sum _{n\geq
0}e_{n}(x_{1},x_{2},\ldots )y^{n}. $$
 
\begin{prop}
\label{elem}
For $k$ nonnegative and $n$ positive,
\begin{eqnarray*}
s(n,k) &=& (-1)^{n}(n-1)!e_{k-1}(-1,-1/2,\ldots ,-1/n)\\*
&=& e_{n-k}(-1,-2,-3,\ldots ,-n+1). \Box  
\end{eqnarray*}
\end{prop}

\digress{We {\em digress} to discuss the implications of
Porism~\ref{complete} and Proposition~\ref{elem}.
The Stirling numbers of the first type with nonnegative degree are merely
examples of the complete 
elementary symmetric function, and those with
 negative degree are merely examples of the elementary
elementary symmetric function. This duality is not too
surprising in light of \cite{MNR} which interprets the
elementary symmetric function as an extension of the
complete symmetric function to a ``negative'' number of
variables. In fact, if we adopt the notation of \cite{hybrid}, then we deduce 
$$  s(n,k) =\lim_{\epsilon \rightarrow 0} \frac{h_{k}
\left(-\frac{1}{1+\epsilon}, -\frac{1}{2+\epsilon }, \ldots
,-\frac{1}{n+\epsilon} \right)}{\Gamma (n+1+\epsilon )}$$ 
 for {\em all} integers $n$ and nonnegative integers $k$.}


\begin{thebibliography}{99}
\bibitem{Char}
{\sc S. Charalambides}, {\em A Review of Stirling Numbers:
Their Generalization and Statistical Application},
Communications in Statistical Theory: Theory and Methods
{\bf 17(A)} 2533-95 (1988).
\bibitem{ch3}
{\sc D. Loeb}, {\em A Generalization of the Binomial Coefficients}, To appear.
\bibitem{thesis}
{\sc D. Loeb}, {\em The Iterated Logarithmic Algebra}, MIT Department of
Mathematics Thesis (1989).
\bibitem{ch4}
{\sc D. Loeb}, {\em The Iterated Logarithmic Algebra}, To appear.
\bibitem{ch5}
{\sc D. Loeb}, {\em The Iterated Logarithmic Algebra II: Sheffer Sequences}, To appear.
\bibitem{hybrid}
{\sc D. Loeb}, {\em Sets with a Negative Number of Elements}, unpublished.
\bibitem{ch6}
{\sc D. Loeb}, {\em Sequences of Symmetric Functions of Binomial Type}, To appear.
\bibitem{LR}
{\sc D. Loeb and G.-C. Rota}, {\em Formal Power Series of
Logarithmic Type}, Advances in Mathematics, {\bf 75} (1989), 1--118.
\bibitem{ch1}
{\sc D. Loeb}, {\em Series with General Exponents}, To appear ???.
\bibitem{MNR}
{\sc N. Metropolis, G. Nicoletti, G.-C. Rota,} 
{\em A New Class of Symmetric Functions,}
Mathematical Analysis and Applications, Part {\bf V}, Volume
{\bf 7B} (1981), 563-575.
\bibitem{Ueno}
{\sc K. Ueno},
{\em Umbral Calculus and Special Functions}, Advances in
Mathematics,
{\bf 67} (1988), 174--229.
\bibitem{UP}
{\sc K. Ueno}, personal correspondence.
\end{thebibliography}
\end{document}